\documentclass[11pt,a4paper]{article}
\usepackage{a4wide}
\usepackage{theorem,paralist}
\usepackage{amsmath,amssymb}
\usepackage[english]{babel}
\usepackage{graphicx}
\usepackage[T1]{fontenc}
\usepackage[numbers]{natbib}
\usepackage{ifpdf}
\usepackage{hyperref}
\usepackage{setspace}

\hypersetup{
    bookmarksopen=false,
    bookmarksnumbered=true,
}
\ifpdf
  \hypersetup{colorlinks=true,linkcolor=black,urlcolor=black,citecolor=black,pdfpagemode=UseOutlines,plainpages=false,pdfpagelabels}
\else
  \hypersetup{colorlinks=true,linkcolor=black,urlcolor=black,citecolor=black,pdfpagemode=UseOutlines,plainpages=false,pdfpagelabels}
\fi

\theorembodyfont{\upshape}
\theoremheaderfont{\bfseries}

\newtheorem{theorem}{Theorem}[section]

\newtheorem{lemma}[theorem]{Lemma}
\newtheorem{remark}[theorem]{Remark}
\newtheorem{proposition}[theorem]{Proposition}

\newenvironment{proof}{\par\textbf{Proof:}\\}{\hfill$\square$\par}
\numberwithin{equation}{section}

\newcommand{\E}{\mathbb{E}}
\renewcommand{\P}{\mathbb{P}}
\renewcommand{\O}{\mathcal{O}}

\newcommand{\PGF}[1]{\widetilde{#1}}
\newcommand{\LST}[1]{\widetilde{#1}}
\newcommand{\dphi}[1]{\frac{\partial\phi_{#1}(\xi,\delta)}{\partial\xi}}
\newcommand{\ddphi}[1]{\frac{\partial^2\phi_{#1}(\xi,\delta)}{\partial\xi^2}}
\newcommand{\taylorphi}[2]{\phi_{#1}(\xi,\delta)#2\delta\dphi{#1}+\frac{\delta^2}{2}\ddphi{#1}}
\providecommand{\href}[2]{#2}

\title{Heavy-traffic analysis of $k$-limited polling systems}
\author{M.A.A. Boon\footnote{Eurandom and Department of Mathematics and Computer Science, Eindhoven University of Technology, P.O. Box 513, 5600MB Eindhoven, The Netherlands}\\\href{mailto:marko@win.tue.nl}{marko@win.tue.nl} \and E.M.M. Winands \footnote{
University of Amsterdam, Korteweg-de Vries Institute for Mathematics, Science Park 904, 1098 XH  Amsterdam, The Netherlands}\\\href{mailto:e.m.m.winands@uva.nl}{e.m.m.winands@uva.nl}}

\date{May, 2014}
\begin{document}
\maketitle

\begin{abstract}
In this paper we study a two-queue polling model with zero switch-over times and $k$-limited service (serve at most $k_i$ customers during one visit period to queue $i$, $i=1,2$) in each queue. The arrival processes at the two queues are Poisson, and the service times are exponentially distributed. By increasing the arrival intensities until one of the queues becomes critically loaded, we derive exact heavy-traffic limits for the joint queue-length distribution using a singular-perturbation technique. It turns out that the number of customers in the stable queue has the same distribution as the number of customers in a vacation system with Erlang-$k_2$ distributed vacations. The queue-length distribution of the critically loaded queue, after applying an appropriate scaling, is exponentially distributed. Finally, we show that the two queue-length processes are independent in heavy traffic.


\end{abstract}

\section{Introduction}\label{sect:introduction}

This paper considers a two-queue $k$-limited  polling model with exponentially distributed service times and zero switch-over times. Under the $k$-limited strategy the server continues working  until either a predefined number of $k_i$ customers is served at queue $i$ or until the queue becomes empty, whichever occurs first.  The interest for this model is fueled by a number of different application areas. That is, the $k$-limited strategy did not only prove its merit in communication systems (see, e.g.,
\cite{borst,charzinski}), but also in the field of logistics (see, e.g. \cite{winands}), and vehicle actuated traffic signals \cite{boontrafficlights2012}. In the present paper, we consider the heavy-traffic scenario, in which one of the queues becomes critically loaded with the other queue remaining stable.

Although the number of papers on polling systems is impressive, hardly any exact results for polling systems with the $k$-limited service policy have been obtained. This can be explained by the fact that the $k$-limited strategy does not satisfy a well-known branching property for polling systems, independently discovered by Fuhrmann \cite{fuhrmann} and Resing \cite{resing}, which significantly increases the analytical complexity. For this reason, most of the papers in the polling literature focus on branching-type service disciplines, like exhaustive service (serve all customers in a queue until it is empty, before switching to the next queue) or gated service (serve all customers present at the server's arrival to a queue, before switching to the next queue).
Groenendijk \cite{Groenendijk1} and Ibe \cite{ibe} give an explicit Laplace-Stieltjes Transform for the waiting-time distribution in a two-queue $1$-limited/exhaustive system. For two-queue systems where both queues are served according to the $1$-limited discipline, the problem of finding the queue length distribution can be shown to translate into a boundary value problem \cite{boxma1,boxma2,cohen1,eisenberg1}.  For general $k$, an exact evaluation for the queue-length distribution is known in two-queue exhaustive/$k$-limited systems (see \cite{lee,ozawa1,ozawa2,winands}). The situation where both queues follow  the $k$-limited discipline has not been solved yet.

When studying the literature on heavy-traffic results for polling systems, a remarkable observation is that distributional heavy-traffic limits have only been rigorously proven for systems with branching-type service disciplines (there is only exception as described later). Even for these systems, the proof is limited to exhaustive systems consisting of two queues \cite{coffman95,coffman98} or systems with Poisson arrivals \cite{RvdM_QUESTA}; all other results are based on conjectures.  That is, a proof for the general setting is already a challenge for almost $20$ years since the publication of the Coffman/Puhalskii/Reiman papers \cite{coffman95,coffman98}. Moreover, the only paper analysing the heavy-traffic behaviour for the $k$-limited discipline is the paper of Lee \cite{lee} for the two-queue exhaustive/$k$-limited system without setup times. More specifically, he studies the limiting regime where the exhaustive queue remains stable, while the $k$-limited queue becomes critically loaded in the limit. Lee \cite{lee} uses results for this system under stable conditions as a basis to derive heavy-traffic asymptotics by exploiting the fact that significant simplifications result as the load increases. However, corresponding steady-state results are not available for the $k$-limited system under consideration in the current paper, which impels us to look for a different family of techniques.

The main contribution of the present paper is that we derive heavy-traffic asymptotics for $k$-limited polling models via the singular-perturbation technique. That is, by increasing the arrival intensities until one of the queues becomes critically loaded, we derive heavy-traffic limits for the joint queue-length distribution in a two-queue $k$-limited  polling model. In this way, we derive the lowest-order asymptotic  to the joint  queue-length distribution in terms of a small positive parameter measuring the closeness of the system to instability.

Furthermore, the results obtained in the present paper provide new insights into the heavy-traffic behaviour of $k$-limited polling systems. It is shown that the number of customers in the stable queue has the same distribution as the number of customers in a vacation system with Erlang-$k_2$ distributed vacations, while the scaled queue-length distribution of the critically loaded queue is exponentially distributed. Finally, we prove that the two queue-length processes are independent in heavy traffic. The singular-perturbation technique can also be extended to an $N$-queue system ($N \geq 2$) with one queue becoming critically loaded. In this limiting regime the stable queues have the same joint distribution of a $k$-limited polling model with $N-1$ queues and an extended switch-over time. These results do not only generalise those derived in \cite{lee} (where the special 2-queue case is studied in which only one queue is served according to the $k$-limited service policy), but are also obtained via a fundamentally different singular-perturbation approach. We would like to note that the observed heavy-traffic behavior is different from the behaviour seen for branching policies in which the cycle times normally converges to a Gamma distribution and the individual workloads follow from an averaging principle.

Finally, to our opinion, the novel application of the singular perturbation technique to polling systems is interesting in itself. That is, the merits of this technique are in its intrinsic simplicity and its intuitively appealing derivation, although it requires some distributional assumptions. The singular perturbation framework is known as a powerful tool to determine asymptotics in all kinds of queueing models, but it has never been applied to polling systems. We would like to refer to Knessl and Tier \cite{knessltiersurvey} for an excellent survey of applications of the perturbation technique to queueing models. It is noteworthy that our paper is inspired by the manner in which Morrison and Borst \cite{morrisonborst2010} apply this technique to a model with interacting queues.

The paper is structured as follows. In the next section we introduce the model and notation. In Section~\ref{sect:analysis} we apply a perturbation technique to study the system under heavy-traffic conditions and derive the limiting scaled joint queue-length distribution. In Section~\ref{sect:conclusions} we interpret the results and give some suggestions on further research. The appendices contain some lengthy derivations required for the analysis in Section~\ref{sect:analysis}.

\section{Model description and notation}\label{sect:model}

We consider a polling model consisting of two queues, $Q_1$ and $Q_2$, that are alternately visited by a single server. Throughout this paper, the subscript~$i$ will always be used to refer to one of the queues, meaning that it always takes on the values 1 or~2. When a server arrives at $Q_i$, it serves at most $k_i$ customers. When $k_i$ customers have been served or $Q_i$ becomes empty, whichever occurs first, the server switches to the other queue. We assume that switching from one queue to the other requires no time. If the other queue turns out to be empty, the server switches back and serves, again, at most $k_i$ customers. If both queues are empty, the server waits until the first arrival and switches to the corresponding queue (say, $Q_j$) to start another visit period of at most $k_j$ customers, $j=1,2$. Customers arrive at $Q_i$ according to a Poisson process with intensity $\lambda_i$. We assume that the service times of customers in $Q_i$ are independent and exponentially distributed with parameter $\mu_i$. We denote the load of the system by $\rho = \rho_1+\rho_2$, where $\rho_i=\lambda_i/\mu_i$. For a polling model without switch-over times, the stability condition is $\rho<1$ \cite{down,fricker}. Furthermore, we assume that
\begin{equation}
\frac{\lambda_1}{k_1} < \frac{\lambda_2}{k_2}.\label{stability2}
\end{equation}
This assumption is discussed in more detail in the next section.

The number of customers in $Q_i$ at time $t$, $t \geq 0$, is denoted by $N_i(t)$. In order to describe the queue length process as a Markov process, we use the approach of Blanc \cite{blanc92}, introducing a supplementary variable $H(t)$ which takes on values $1,2,\dots,k_1+k_2$. The variable $H(t)$ is used to determine the server position ($Q_1$ or $Q_2$) at time $t$ and the number of customers that can be served before the server has to switch to the next queue. When $1\leq H(t) \leq k_1$, this means that the server is serving $Q_1$ at time $t$, and that the customer in service is the $H(t)$-th customer being served during the present visit period. If $k_1+1 \leq H(t) \leq k_1+k_2$, the server is serving the ($H(t)-k_1$)-th customer in $Q_2$. Now, $(N_1(t), N_2(t), H(t))$ is a Markov process. Assuming that $\rho<1$, define the stationary probabilities $p(n_1, n_2, h) := \lim_{t\rightarrow\infty}\P\big(N_1(t)=n_1,N_2(t)=n_2,H(t)=h\big)$, and define the steady-state queue lengths $N_i$.

\section{Analysis}\label{sect:analysis}

We study the heavy-traffic limit of the joint queue-length process $(N_1, N_2)$ by increasing the arrival rate $\lambda_2$, while keeping $\lambda_1$ fixed. When $\rho$ tends to 1, Assumption \eqref{stability2} implies that $Q_2$ will become critically loaded, whereas $Q_1$ remains stable due to the fact that at most $k_2$ customers are served during each visit period at $Q_2$. In case $\lambda_1/k_1=\lambda_2/k_2$ in the limit, both queues would become critically loaded simultaneously and the system behaviour is different from the limiting behaviour found in the present paper. In fact, the limiting queue-length behaviour for that specific case remains an open problem. We discuss this topic briefly in Section \ref{sect:conclusions}.

We use a single-perturbation method to find the queue-length distributions in heavy-traffic $(\rho \uparrow 1)$. First, we write down the balance equations of our model and apply a perturbation to the arrival rate of $Q_2$ to these equations, in the case that this queue is close to becoming critically loaded. By solving the system of balance equations for the lowest (zeroth) order terms, we find the queue length distribution of the stable queue, $Q_1$. By analysing the equations for the first-order and second-order terms, we also obtain a differential equation for the scaled number of customers in $Q_2$, which we can solve to show that this number converges to an exponential distribution.

\subsection{Balance equations}

The balance equations for a polling model with exponentially distributed service times and $k$-limited service at each of the queues are given by Blanc \cite{blanc92}. For completeness, we present these equations for our two-queue model below.

\begin{subequations}
\begin{alignat}{3}
&(\lambda_1+\lambda_2+\mu_2)p(0,n_2,k_1+1) &&= \lambda_2p(0,n_2-1,k_1+1)+\mu_2p(0, n_2+1,k_1+k_2)\nonumber\\
&                                          &&+\sum_{h=1}^{k_1} \mu_1p(1,n_2,h),\label{balanceqn1}\\
&(\lambda_1+\lambda_2+\mu_2)p(0,n_2, h_2) &&= \lambda_2p(0,n_2-1,h_2)+\mu_2p(0,n_2+1,h_2-1), \\
&(\lambda_1+\lambda_2+\mu_1)p(1,n_2,1) &&= \lambda_2p(1,n_2-1,1)+\mu_2p(1,n_2+1,k_1+k_2),\\
&(\lambda_1+\lambda_2+\mu_1)p(1,n_2,h_1) &&= \lambda_2p(1,n_2-1,h_1)+\mu_1p(2,n_2,h_1-1), \\
&(\lambda_1+\lambda_2+\mu_1)p(n_1+1,n_2,1) &&= \lambda_1p(n_1,n_2,1)+\lambda_2p(n_1+1,n_2-1,1)\nonumber\\
&                                        &&+\mu_2p(n_1+1,n_2+1,k_1+k_2),\\
&(\lambda_1+\lambda_2+\mu_1)p(n_1+1,n_2,h_1) &&= \lambda_1p(n_1,n_2,h_1)+\lambda_2p(n_1+1,n_2-1,h_1)\nonumber\\
&                                          &&+\mu_1p(n_1+2,n_2,h_1-1),\\
&(\lambda_1+\lambda_2+\mu_2)p(n_1,n_2,k_1+1) &&= \lambda_1p(n_1-1,n_2,k_1+1)+\lambda_2p(n_1,n_2-1,k_1+1)\nonumber\\
&                                            &&+\mu_1p(n_1+1,n_2,k_1),\\
&(\lambda_1+\lambda_2+\mu_2)p(n_1,n_2,h_2) &&= \lambda_1p(n_1-1,n_2,h_2)+\lambda_2p(n_1,n_2-1,h_2)\nonumber\\
&                                          &&+\mu_2p(n_1,n_2+1,h_2-1),\label{balanceqnlast}
\end{alignat}
\end{subequations}
for $n_1=1,2,\dots; n_2=2,3,\dots; h_1=2,3,\dots,k_1$, and $h_2=k_1+2,\dots,k_1+k_2$. Note that \eqref{balanceqn1}-\eqref{balanceqnlast} are \emph{not all} balance equations. We have omitted all equations for $n_2=0$ and $n_2=1$, since it will turn out that these do not play a role after the perturbation. The intuitive explanation is that $N_2(t)$ will tend to infinity as $Q_2$ becomes critically loaded and the probabilities $p(n_1,n_2,h)$ become negligible for low values of $n_2$.

\subsection{Perturbation}

From the stability condition we have that the system becomes unstable as ${\lambda_2}/{\mu_2}\uparrow1-{\lambda_1}/{\mu_1}$, which means that the arrival rate $\lambda_2$ approaches $\mu_2(1-\lambda_1/\mu_1)$. Therefore we will assume that
\begin{equation}
\lambda_2 = \mu_2\left(1-\frac{\lambda_1}{\mu_1}\right)-\delta\omega,\qquad \omega > 0, 0<\delta \ll 1.\label{perturbation}
\end{equation}
At the end of this section we will take an appropriate choice for the constant $\omega$, which will influence the limit of the scaled queue length in $Q_2$.

Let $\xi=\delta n_2$, and
\begin{equation}
p(n_1, \xi/\delta, h) = \delta \phi_{n_1, h}(\xi,\delta), \qquad 0<\xi=\O(1), h=1,2,\dots,k_1+k_2.\label{phi}
\end{equation}
Note that once $\omega$ is chosen, $\delta$ and the scaled variable $\xi$ will be uniquely defined.
The next step is to substitute \eqref{perturbation} and \eqref{phi} in the balance equations \eqref{balanceqn1}-\eqref{balanceqnlast}, and take the Taylor series expansion with respect to $\delta$. For reasons of compactness, we only show the intermediate results for Equation \eqref{balanceqnlast} as an illustration:
\begin{multline*}
(\lambda_1+\mu_2)\phi_{n_1,h_2}(\xi,\delta) - \lambda_1\phi_{n_1-1,h_2}(\xi,\delta)-\mu_2\phi_{n_1,h_2-1}(\xi+\delta,\delta)=\\
\Big(\mu_2\big(1-\frac{\lambda_1}{\mu_1}\big)-\delta\omega\Big)\big(\phi_{n_1,h_2}(\xi-\delta,\delta)-\phi_{n_1,h_2}(\xi,\delta)\big).
\end{multline*}
Taking the Taylor series yields:
\begin{multline}
(\lambda_1+\mu_2)\phi_{n_1,h_2}(\xi,\delta) - \lambda_1\phi_{n_1-1,h_2}(\xi,\delta)\\
-\mu_2\left(\taylorphi{n_1,h_2-1}{+}\right)=\\
-\Big(\mu_2\big(1-\frac{\lambda_1}{\mu_1}\big)-\delta\omega\Big)\left(\delta\dphi{n_1,h_2}-\frac{\delta^2}{2}\ddphi{n_1,h_2}\right) +\O(\delta^3).
\label{taylor}
\end{multline}
Note that $\lambda_2$ (or $\mu_2\big(1-\frac{\lambda_1}{\mu_1}\big)-\delta\omega$ after the substitution) only plays a role in this equation for $\O(\delta)$ terms and higher. It is readily verified that this is the case for all balance equations. We now expand in powers of $\delta$, and let
\begin{equation}
\phi_{n_1,h}(\xi,\delta) = \phi^{(0)}_{n_1,h}(\xi) +\delta\phi^{(1)}_{n_1,h}(\xi) +\O(\delta^2).\label{phi0subst}
\end{equation}
One would expect that also $\delta^2\phi^{(2)}_{n_1,h}(\xi)$ should be introduced in order to analyse the second-order terms. However, Proposition \ref{prop1} implies that the above power expansion suffices.
We also define the corresponding generating functions
\begin{align*}
\PGF{Q}_h(z, \xi, \delta)&:=\sum_{n_1=0}^\infty\phi_{n_1,h}(\xi,\delta)z^{n_1}, &&\quad \PGF{Q}^{(j)}_h(z, \xi):=\sum_{n_1=0}^\infty\phi_{n_1,h}^{(j)}(\xi)z^{n_1}, &&\quad j=0,1,2,\dots.
\end{align*}

In the next subsections we first equate the lowest order terms of the resulting equations to find an expression for (the generating function of) $\phi^{(0)}_{n,h}(\xi)$, and subsequently we equate the first-order and second-order terms to find the scaled queue-length distribution of $Q_2$.

\subsection{Equating the lowest-order terms}\label{sect:equatingO1}

Equating the lowest-order terms of the balance equations, after substituting \eqref{perturbation}, \eqref{phi}, and \eqref{phi0subst}, results in the following equations.
\begin{subequations}
\begin{alignat}{3}
&(\lambda_1+\mu_2)\phi^{(0)}_{0,k_1+1}(\xi) &&= \mu_2 \phi^{(0)}_{0,k_1+k_2}(\xi)+\sum_{h=1}^{k_1} \mu_1\phi^{(0)}_{1,h}(\xi),\label{balanceqn1order0}\\
&(\lambda_1+\mu_2)\phi^{(0)}_{0,h_2}(\xi) &&= \mu_2\phi^{(0)}_{0,h_2-1}(\xi),\label{balanceqn2order0} \\
&(\lambda_1+\mu_1)\phi^{(0)}_{1,1}(\xi) &&= \mu_2\phi^{(0)}_{1,k_1+k_2}(\xi),\label{balanceqn3order0}\\
&(\lambda_1+\mu_1)\phi^{(0)}_{1,h_1}(\xi) &&= \mu_1\phi^{(0)}_{2,h_1-1}(\xi),\label{balanceqn4order0} \\
&(\lambda_1+\mu_1)\phi^{(0)}_{n_1+1,1}(\xi) &&= \lambda_1\phi^{(0)}_{n_1,1}(\xi)+\mu_2\phi^{(0)}_{n_1+1,k_1+k_2}(\xi),\label{balanceqn5order0}\\
&(\lambda_1+\mu_1)\phi^{(0)}_{n_1+1,h_1}(\xi) &&= \lambda_1\phi^{(0)}_{n_1,h_1}(\xi)+\mu_1\phi^{(0)}_{n_1+2,h_1-1}(\xi),\label{balanceqn6order0}\\
&(\lambda_1+\mu_2)\phi^{(0)}_{n_1,k_1+1}(\xi) &&= \lambda_1\phi^{(0)}_{n_1-1,k_1+1}(\xi)+\mu_1\phi^{(0)}_{n_1+1,k_1}(\xi),\label{balanceqn7order0}\\
&(\lambda_1+\mu_2)\phi^{(0)}_{n_1,h_2}(\xi) &&= \lambda_1\phi^{(0)}_{n_1-1,h_2}(\xi)+\mu_2\phi^{(0)}_{n_1,h_2-1}(\xi),\label{balanceqnlastorder0}
\end{alignat}
\end{subequations}
for $n_1=1,2,\dots; h_1=2,3,\dots,k_1$, and $h_2=k_1+2,\dots,k_1+k_2$.

Note that $\sum_{n_1=0}^\infty\sum_{h=1}^{k_1+k_2}\phi^{(0)}_{n_1,h}(\xi) \neq 1$. For this reason we introduce $P_0(\xi)$ and $\pi^{(0)}_{n_1,h}$, with
\begin{equation}
\phi^{(0)}_{n_1,h}(\xi) = \pi^{(0)}_{n_1,h}P_0(\xi), \qquad \textrm{ and } \qquad \sum_{n_1=0}^\infty\sum_{h=1}^{k_1+k_2}\pi^{(0)}_{n_1,h} = 1,\label{phi0pi0subst}
\end{equation}
for $n_1=0,1,2,\dots$ and $h=1,2,\dots,k_1+k_2$.
Careful inspection of these balance equations reveals that equations \eqref{balanceqn1order0}-\eqref{balanceqnlastorder0} describe the behaviour of a single-server vacation queue with the following properties:
\begin{compactenum}[P1.]
\item the arrival process is Poisson with intensity $\lambda_1$,\label{vacationprop1}
\item the service times are exponentially distributed with mean $1/\mu_1$,
\item the service discipline is $k$-limited service with service limit $k_1$,
\item the vacations are Erlang-$k_2$ distributed with parameter $\mu_2$,
\item whenever the server finds the system empty upon return from a vacation, it immediately starts another vacation.\label{vacationproplast}
\end{compactenum}
This system has been studied in the literature (cf. \cite{lee89}) and, in general, no closed-form expressions for the steady-state queue-length probabilities can be obtained. However, it is possible to find the probability generating function (PGF) of the queue-length distribution. Define
\[
\LST{L}^{(0)}(z) := \sum_{h=1}^{k_1+k_2} \LST{L}^{(0)}_h(z), \qquad\textrm{ where }\qquad \LST{L}^{(0)}_h(z) := \sum_{n=0}^\infty \pi^{(0)}_{n,h}z^{n}, \qquad h=1,\dots,k_1+k_2,
\]
and
\begin{equation}
\LST{G}(z) := \frac{\mu_1}{\lambda_1(1-z)+\mu_1}, \qquad
\qquad
\LST{H}(z) := \frac{\mu_2}{\lambda_1(1-z)+\mu_2}.\label{ghz}
\end{equation}
It is easily seen that $\LST{G}(z)$ and $\LST{H}(z)$ are the PGFs of the number of arrivals during respectively one service, and during one stage of the vacation (which consists of $k_2$ exponential stages). It follows from \eqref{balanceqn1order0}-\eqref{balanceqnlastorder0} that
\begin{align}
\LST{L}^{(0)}_{k_1+k_2}(z) &= \frac{\displaystyle
\LST{H}(z)^{k_2}\left[\pi^{(0)}_{0,k_1+k_2} \left(1- \left({\LST{G}(z)}/{z}\right)^{k_1}\right) +\frac{\mu_1}{\mu_2}\sum_{h=1}^{k_1-1}\pi_{1,h}^{(0)}\left(1-\left({\LST{G}(z)}/{z}\right)^{k_1-h} \right) \right]}{1-\big(\LST{G}(z)/z\big)^{k_1}\LST{H}(z)^{k_2}},\label{lstL1}\\
\LST{L}^{(0)}_{h_1}(z) &= \frac{z\mu_2}{\mu_1}\big(\LST{G}(z)/z\big)^{h_1}\left[\LST{L}^{(0)}_{k_1+k_2}(z)-\pi^{(0)}_{0,k_1+k_2}\right]-\sum_{h=1}^{h_1-1}\pi_{1,h}^{(0)}z\left({\LST{G}(z)}/{z}\right)^{h_1-h}, \label{lstL2}\\
\LST{L}^{(0)}_{h_2}(z) &= \LST{H}(z)^{h_2-k_1}\left[\frac{\mu_1}{\mu_2}\sum_{h=1}^{k_1-1}\pi^{(0)}_{1,h}+\pi^{(0)}_{0,k_1+k_2}      +\frac{\mu_1}{z\mu_2}\LST{L}^{(0)}_{k_1}(z)\right], \label{lstL3}
\end{align}
for $h_1=1,\dots,k_1$ and $h_2=k_1+1,\dots,k_1+k_2-1$. A derivation of \eqref{lstL1}-\eqref{lstL3} can be found in Appendix \ref{proofVacation}.
These equations still contain $k_1$ unknowns: $\pi^{(0)}_{1,1},\pi^{(0)}_{1,2},\dots,\pi^{(0)}_{1,k_1-1},$ and $\pi^{(0)}_{0,k_1+k_2}$. See Appendix \ref{proofVacation} for more details on how to eliminate them using Rouch\'e's Theorem. Foregoing the derivation of the limiting behaviour of $Q_2$ we already would like to mention that these unknowns do not play a role therein.

We conclude from equating the lowest-order terms of the balance equations \eqref{balanceqn1}-\eqref{balanceqnlast}, after substituting \eqref{perturbation}, \eqref{phi}, and \eqref{phi0subst}, that
\begin{equation}
\sum_{n_1=0}^\infty\sum_{h=1}^{k_1+k_2}\phi^{(0)}_{n_1,h}(\xi)z^{n_1} = \LST{L}^{(0)}(z)P_0(\xi) ,\label{solutionOrder0}
\end{equation}
where $P_0(\xi)$ still has to be determined. Consequently, in heavy traffic, the queue length of the stable queue ($Q_1$) has the same distribution as the queue length in a vacation system with Erlang($k_2$) distributed vacations with parameter $\mu_2$, exponential service times with parameter $\mu_1$, and $k_1$-limited service. The assumption that we have Poisson arrivals and first-come-first-served service implies that we can use the distributional form of Little's law to obtain (the Laplace-Stieltjes transform of) the waiting-time distribution of customers in $Q_1$ (see, for example, Keilson and Servi \cite{keilsonservi90}).

\begin{remark}
In this paper we assume that $\lambda_1/k_1 < \lambda_2/k_2$, causing $Q_2$ to become critically loaded when $\lambda_2$ is being increased. We have implicitly used this assumption when solving the balance equations \eqref{balanceqn1}-\eqref{balanceqnlast}. It is well-known that the vacation system described by these equations is stable if and only if
\begin{equation}
\lambda_1\E[C] < k_1,
\label{stabilityvacation}
\end{equation}
where $\E[C]$ is the mean cycle time, i.e., the mean length of one visit period plus one vacation. Denoting the length of a vacation by $S$, we have $\E[C] = \E[S]/(1-\rho_1) = k_2/(\mu_2(1-\lambda_1/\mu_1))$. When substituting this in \eqref{stabilityvacation}, we indeed obtain exactly the same inequality as \eqref{stability2} after substituting \eqref{perturbation} and letting $\delta\downarrow0$.
\end{remark}

\begin{remark}\label{lambda1remark}
Another interesting observation is that one could consider more general ways of varying the arrival rates in order to let $Q_2$ become critically loaded. To this end, we introduce $\lambda_1^*$ and $\lambda_2^*$ such that
$\lambda_1^*/\mu_1+\lambda_2^*/\mu_2 = 1$. Additionally, we assume that
\begin{equation}
\frac{\lambda_1^*}{k_1}<\frac{1}{\frac{k_1}{\mu_1}+\frac{k_2}{\mu_2}}, \textrm{ or equivalently: }\frac{\lambda_2^*}{k_2}>\frac{1}{\frac{k_1}{\mu_1}+\frac{k_2}{\mu_2}}.
\label{stability2alt}
\end{equation}
We now let $\lambda_1\rightarrow\lambda_1^*$ and $\lambda_2\rightarrow\lambda_2^*$ for $\delta\downarrow 0$, with
\begin{equation}
\frac{\lambda_1}{\mu_1} + \frac{\lambda_2}{\mu_2} = 1-\delta\omega^*,\qquad \omega^* > 0, 0<\delta \ll 1.\label{perturbationalt}
\end{equation}
Any arbitrary way in which we let $\lambda_1$ and $\lambda_2$ approach respectively $\lambda_1^*$ and $\lambda_2^*$, for $\delta\downarrow 0$, will cause $Q_2$ to become critically loaded (because of assumption \eqref{stability2alt}). All results obtained in this paper will still be valid, by choosing $\omega^* = \omega/\mu_2$.
\end{remark}

\subsection{Equating the first-order terms}\label{sect:equatingOdelta}

In this section we study, and solve, the system of equations that results from equating the first-order terms of the perturbed balance equations. For notational reasons, we define
\[
\psi^{(1)}_{n_1,h}(\xi) := \phi^{(1)}_{n_1,h}(\xi) + \phi'^{\,(0)}_{n_1,h}(\xi),\qquad\textrm{ where }
\qquad\phi'^{\,(0)}_{n_1,h}(\xi) := \frac{d \phi^{(0)}_{n_1,h}(\xi)}{d \xi},
\]
for $n_1=0,1,\dots$ and $h=1,2,\dots,k_1+k_2$. The resulting set of equations for the probabilities $\phi^{(1)}_{n_1,h}(\xi)$ is given below.

\begin{subequations}
\begin{alignat}{3}
&(\lambda_1+\mu_2)\phi^{(1)}_{0,k_1+1}(\xi) &&= \mu_2 \psi^{(1)}_{0,k_1+k_2}(\xi)+\sum_{h=1}^{k_1} \mu_1\phi^{(1)}_{1,h}(\xi) -\mu_2\left(1-\frac{\lambda_1}{\mu_1}\right)\phi'^{\,(0)}_{0,k_1+1}(\xi),\label{balanceqn1order1}\\
&(\lambda_1+\mu_2)\phi^{(1)}_{0,h_2}(\xi) &&= \mu_2\psi^{(1)}_{0,h_2-1}(\xi) -\mu_2\left(1-\frac{\lambda_1}{\mu_1}\right)\phi'^{\,(0)}_{0,h_2}(\xi),\label{balanceqn2order1} \\
&(\lambda_1+\mu_1)\phi^{(1)}_{1,1}(\xi) &&= \mu_2\psi^{(1)}_{1,k_1+k_2}(\xi) -\mu_2\left(1-\frac{\lambda_1}{\mu_1}\right)\phi'^{\,(0)}_{1,1}(\xi),\label{balanceqn3order1}\\
&(\lambda_1+\mu_1)\phi^{(1)}_{1,h_1}(\xi) &&= \mu_1\phi^{(1)}_{2,h_1-1}(\xi) -\mu_2\left(1-\frac{\lambda_1}{\mu_1}\right)\phi'^{\,(0)}_{1,h_1}(\xi),\label{balanceqn4order1} \\
&(\lambda_1+\mu_1)\phi^{(1)}_{n_1+1,1}(\xi) &&= \lambda_1\phi^{(1)}_{n_1,1}(\xi)+\mu_2\psi^{(1)}_{n_1+1,k_1+k_2}(\xi) -\mu_2\left(1-\frac{\lambda_1}{\mu_1}\right)\phi'^{\,(0)}_{n_1+1,1}(\xi),\label{balanceqn5order1}\\
&(\lambda_1+\mu_1)\phi^{(1)}_{n_1+1,h_1}(\xi) &&= \lambda_1\phi^{(1)}_{n_1,h_1}(\xi)+\mu_1\phi^{(1)}_{n_1+2,h_1-1}(\xi) -\mu_2\left(1-\frac{\lambda_1}{\mu_1}\right)\phi'^{\,(0)}_{n_1+1,h_1}(\xi),\label{balanceqn6order1}\\
&(\lambda_1+\mu_2)\phi^{(1)}_{n_1,k_1+1}(\xi) &&= \lambda_1\phi^{(1)}_{n_1-1,k_1+1}(\xi)+\mu_1\phi^{(1)}_{n_1+1,k_1}(\xi) -\mu_2\left(1-\frac{\lambda_1}{\mu_1}\right)\phi'^{\,(0)}_{n_1,k_1+1}(\xi),\label{balanceqn7order1}\\
&(\lambda_1+\mu_2)\phi^{(1)}_{n_1,h_2}(\xi) &&= \lambda_1\phi^{(1)}_{n_1-1,h_2}(\xi)+\mu_2\psi^{(1)}_{n_1,h_2-1}(\xi) -\mu_2\left(1-\frac{\lambda_1}{\mu_1}\right)\phi'^{\,(0)}_{n_1,h_2}(\xi),\label{balanceqnlastorder1}
\end{alignat}
\end{subequations}
for $n_1=1,2,\dots; h_1=2,3,\dots,k_1$, and $h_2=k_1+2,\dots,k_1+k_2$. The solution to this system of equations, in terms of generating functions, can be found in Appendix \ref{proofOrder2}. This solution is used to derive the following relation
\begin{equation}
\frac{\lambda_1}{\mu_1}\sum_{n=0}^\infty\sum_{h=1}^{k_1+k_2}\phi^{(1)}_{n,h}(\xi) - \sum_{n=0}^\infty\sum_{h_1=1}^{k_1}\phi^{(1)}_{n,h_1}(\xi)
= \frac{\lambda_1\mu_2}{\mu_1^2}P_0'(\xi),\label{relation2ndOrder}
\end{equation}
which turns out to play a key role in determining the HT limit of the joint queue-length distribution (see the next section).

\subsection{Equating the second-order terms}\label{sect:equatingOdelta2}

In order to find an expression for $P_0(\xi)$ and, consequently, solve \eqref{solutionOrder0}, we consider the Taylor series of all perturbed balance equations. In this section we show that, fortunately, we only need to consider the \emph{sum} of all these equations (such as Equation \eqref{taylor}) over all $n_1=0,1,2,\dots$ and $h=1,2,\dots,k_1+k_2$, and we consecutively consider the $\O(1)$, $\O(\delta)$, and $\O(\delta^2)$ terms. Using the results we have obtained so far, we prove that:
\begin{enumerate}
\item all $\O(1)$, i.e., all $\phi_{n_1,h}(\xi,\delta)$, terms cancel immediately,
\item the $\O(\delta)$ terms cancel after expanding $\phi_{n_1,h}(\xi,\delta)$ in powers of $\delta$ (i.e., substituting \eqref{phi0subst}),
\item the equation that results from equating the $\O(\delta^2)$ terms, can be solved to find $P_0(\xi)$.
\end{enumerate}
The above three results are proven in Propositions \ref{prop1}, \ref{prop2}, and \ref{prop3}.

\begin{proposition}\label{prop1}
After taking the summation over all $n_1=0,1,2,\dots$ and $h=1,2,\dots,k_1+k_2$ of the Taylor series of all perturbed balance equations, the $\O(1)$ terms cancel.
\begin{proof}
We can follow the generating function approach, used in Section \ref{sect:equatingO1}, but replace the probabilities $\pi^{(0)}_{n_1,h}$, by $\phi_{n_1,h}(\xi,\delta)$. This results in the same set of equations as \eqref{lstL1}-\eqref{lstL3}, but with terms $\PGF{Q}_h(z, \xi, \delta)$ and $\phi_{n_1,h}(\xi,\delta)$ instead of $\PGF{L}^{(0)}_h(z)$ and $\pi^{(0)}_{n_1,h}$. Substituting $z=1$ yields:
\begin{align}
\mu_2\PGF{Q}_{h_2}(1, \xi, \delta) &= \mu_2\PGF{Q}_{h_2-1}(1, \xi, \delta),\label{qeqn1}\\
\mu_2\PGF{Q}_{k_1+1}(1, \xi, \delta) &= \mu_1\PGF{Q}_{k_1}(1, \xi, \delta)+\mu_1\sum_{h=1}^{k_1-1}\phi_{1,h}(\xi,\delta)+\mu_2\phi_{0,k_1+k_2}(\xi,\delta), \\
\mu_1\PGF{Q}_{h_1}(1, \xi, \delta) &= \mu_1\PGF{Q}_{h_1-1}(1, \xi, \delta)-\mu_1\phi_{1,h_1-1}(\xi,\delta),\\
\mu_1\PGF{Q}_{1}(1, \xi, \delta) &= \mu_2\PGF{Q}_{k_1+k_2}(1, \xi, \delta)-\mu_2\phi_{0,k_1+k_2}(\xi,\delta),\label{qeqnlast}
\end{align}
for $h_1=2,\dots,k_1$ and $h_2=k_1+2,\dots,k_1+k_2$. The summation of \eqref{qeqn1}-\eqref{qeqnlast} over all $h=1,2,\dots,k_1+k_2$ cancels all terms.
\end{proof}
\end{proposition}

\begin{proposition}\label{prop2}
After taking the summation over all $n_1=0,1,2,\dots$ and $h=1,2,\dots,k_1+k_2$ of the Taylor series of all perturbed balance equations, substituting \eqref{phi0subst} and using the results from Appendix \ref{proofVacation}, the $\O(\delta)$ terms cancel.
\begin{proof}
Define
\[
\phi'_{n_1,h}(\xi,\delta) := \frac{\partial \phi_{n_1,h}(\xi,\delta)}{\partial \xi},
\]
for $n_1=0,1,\dots$ and $h=1,2,\dots,k_1+k_2$. Given the fact that the $\O(1)$ terms cancel, taking the $\O(\delta)$ terms leads to the following equations:
\begin{align*}
0 &= \mu_2\phi'_{0,k_1+k_2}(\xi,\delta) - \mu_2\left(1-\frac{\lambda_1}{\mu_1}\right)\phi'_{0,k_1+1}(\xi,\delta),\\
0 &= \mu_2\phi'_{0,h_2-1}(\xi,\delta) - \mu_2\left(1-\frac{\lambda_1}{\mu_1}\right)\phi'_{0,h_2}(\xi,\delta),\\
0 &= \mu_2\phi'_{1,k_1+k_2}(\xi,\delta) - \mu_2\left(1-\frac{\lambda_1}{\mu_1}\right)\phi'_{1,1}(\xi,\delta),\\
0 &= - \mu_2\left(1-\frac{\lambda_1}{\mu_1}\right)\phi'_{1,h_1}(\xi,\delta),\\
0 &= \mu_2\phi'_{n_1+1,k_1+k_2}(\xi,\delta) - \mu_2\left(1-\frac{\lambda_1}{\mu_1}\right)\phi'_{n_1+1,1}(\xi,\delta),\\
0 &= - \mu_2\left(1-\frac{\lambda_1}{\mu_1}\right)\phi'_{n_1+1,h_1}(\xi,\delta),\\
0 &= - \mu_2\left(1-\frac{\lambda_1}{\mu_1}\right)\phi'_{n_1,k_1+1}(\xi,\delta),\\
0 &= \mu_2\phi'_{n_1,h_2-1}(\xi,\delta) - \mu_2\left(1-\frac{\lambda_1}{\mu_1}\right)\phi'_{n_1,h_2}(\xi,\delta),
\end{align*}
for $h_1=2,\dots,k_1$ and $h_2=k_1+2,\dots,k_1+k_2$. Taking the generating functions of these equations and substituting $z=1$ results in the following set of equations:
\begin{align}
0 &= \mu_2\phi'_{0,k_1+k_2}(\xi,\delta)-\mu_2\left(1-\frac{\lambda_1}{\mu_1}\right) \PGF{Q}'_{k_1+1}(1, \xi, \delta),\label{dqeqn1}\\
0 &= \mu_2\PGF{Q}'_{h_2-1}(1, \xi, \delta)-\mu_2\left(1-\frac{\lambda_1}{\mu_1}\right) \PGF{Q}'_{h_2}(1, \xi, \delta),\\
0 &= \mu_2\left(\PGF{Q}'_{k_1+k_2}(1, \xi, \delta)-\phi'_{0,k_1+k_2}(\xi,\delta)\right)-\mu_2\left(1-\frac{\lambda_1}{\mu_1}\right) \PGF{Q}'_{1}(1, \xi, \delta),\\
0 &= -\mu_2\left(1-\frac{\lambda_1}{\mu_1}\right) \PGF{Q}'_{h_1}(1, \xi, \delta),\label{dqeqnlast}
\end{align}
for $h_1=2,\dots,k_1$ and $h_2=k_1+2,\dots,k_1+k_2$. The PGF $\PGF{Q}'_h(z, \xi, \delta)$ is the derivative of $\PGF{Q}_h(z, \xi, \delta)$ with respect to $\xi$.

The summation of \eqref{dqeqn1}-\eqref{dqeqnlast} over all $h=1,2,\dots,k_1+k_2$ yields:
\begin{equation}
\mu_2\sum_{h_1=1}^{k_1}\PGF{Q}'_{h_1}(1, \xi, \delta) = \mu_2\frac{\lambda_1}{\mu_1}\sum_{h=1}^{k_1+k_2}\PGF{Q}'_{h}(1, \xi, \delta).\label{solorder1}
\end{equation}
Apparently, the $\O(\delta)$ terms do not cancel (yet). However, after substituting \eqref{phi0subst}, taking the $\O(\delta)$ terms, and using \eqref{phi0pi0subst}, we obtain:
\begin{equation}
\mu_2\sum_{h_1=1}^{k_1}\PGF{L}^{\,(0)}_{h_1}(1)P_0'(\xi) = \mu_2\frac{\lambda_1}{\mu_1}\sum_{h=1}^{k_1+k_2}\PGF{L}^{\,(0)}_{h}(1)P_0'(\xi).
\label{laccent0}
\end{equation}
Since the (at this moment still unknown) terms $P_0'(\xi)$ cancel out, and since $\sum_{h=1}^{k_1+k_2}\PGF{L}^{\,(0)}_{h}(1)=1$, Equation \eqref{laccent0} reduces to
\[
\sum_{h_1=1}^{k_1}\PGF{L}^{\,(0)}_{h_1}(1) = \rho_1,
\]
which is indeed true (see \eqref{rhovacationsystem}).
\end{proof}
\end{proposition}

\begin{proposition}\label{prop3}
Taking the summation over all $n_1=0,1,2,\dots$ and $h=1,2,\dots,k_1+k_2$ of the Taylor series of all perturbed balance equations and equating the $\O(\delta^2)$ terms, yields the following differential equation for $P_0(\xi)$:
\begin{equation}
\omega P_0'(\xi)=-\left(\mu_2+\frac{\lambda_1\mu_2(\mu_2-\mu_1)}{\mu_1^2}\right)P_0''(\xi).\label{differentialeqn}
\end{equation}
\begin{proof}
Define
\[
\phi''_{n_1,h}(\xi,\delta) := \frac{\partial^2 \phi^{(0)}_{n_1,h}(\xi,\delta)}{\partial \xi^2},
\]
for $n_1=0,1,\dots$ and $h=1,2,\dots,k_1+k_2$. As before, the $\O(1)$ terms cancel. From the proof of the Proposition \ref{prop2} we have learned to include $\O(\delta)$ terms as well, because multiplied by $\delta\phi'^{\,(1)}_{n,h}(\xi)$ these terms are $\O(\delta^2)$ as well. This leads to the following equations:
\begin{align*}
0 =& \mu_2\left[\phi'^{\,(1)}_{0,k_1+k_2}(\xi)+\frac12 \phi''^{\,(0)}_{0,k_1+k_2}(\xi)\right] - \mu_2\left(1-\frac{\lambda_1}{\mu_1}\right)\left[\phi'^{\,(1)}_{0,k_1+1}(\xi)-\frac12\phi''^{\,(0)}_{0,k_1+1}(\xi)\right]
\\
&+\omega\phi'^{\,(0)}_{0,k_1+1}(\xi,\delta),\\
0 =& \mu_2\left[\phi'^{\,(1)}_{0,h_2-1}(\xi)+\frac12 \phi''^{\,(0)}_{0,h_2-1}(\xi)\right]- \mu_2\left(1-\frac{\lambda_1}{\mu_1}\right)\left[\phi'^{\,(1)}_{0,h_2}(\xi)-\frac12\phi''^{\,(0)}_{0,h_2}(\xi)\right] +\omega\phi'^{\,(0)}_{0,h_2}(\xi),\\
0 =& \mu_2\left[\phi'^{\,(1)}_{1,k_1+k_2}(\xi)+\frac12\phi''^{\,(0)}_{1,k_1+k_2}(\xi)\right] - \mu_2\left(1-\frac{\lambda_1}{\mu_1}\right)\left[\phi'^{\,(1)}_{1,1}(\xi)-\frac12\phi''^{\,(0)}_{1,1}(\xi)\right] +\omega\phi'^{\,(0)}_{1,1}(\xi),\\
0 =& -\mu_2\left(1-\frac{\lambda_1}{\mu_1}\right)\left[\phi'^{\,(1)}_{1,h_1}(\xi)-\frac12\phi''^{\,(0)}_{1,h_1}(\xi)\right] +\omega\phi'^{\,(0)}_{1,h_1}(\xi),\\
0 =& \mu_2\left[\phi'^{\,(1)}_{n_1+1,k_1+k_2}(\xi)+\frac12\phi''^{\,(0)}_{n_1+1,k_1+k_2}(\xi)\right] - \mu_2\left(1-\frac{\lambda_1}{\mu_1}\right)\left[\phi'^{\,(1)}_{n_1+1,1}(\xi)-\frac12\phi''^{\,(0)}_{n_1+1,1}(\xi)\right] \\ &+\omega\phi'^{\,(0)}_{n_1+1,1}(\xi),\\
0 =& -\mu_2\left(1-\frac{\lambda_1}{\mu_1}\right)\left[\phi'^{\,(1)}_{n_1+1,h_1}(\xi)-\frac12\phi''^{\,(0)}_{n_1+1,h_1}(\xi)\right] +\omega\phi'^{\,(0)}_{n_1+1,h_1}(\xi),\\
0 =& -\mu_2\left(1-\frac{\lambda_1}{\mu_1}\right)\left[\phi'^{\,(1)}_{n_1,k_1+1}(\xi)-\frac12\phi''^{\,(0)}_{n_1,k_1+1}(\xi)\right] +\omega\phi'^{\,(0)}_{n_1,k_1+1}(\xi),\\
0 =& \mu_2\left[\phi'^{\,(1)}_{n_1,h_2-1}(\xi)+\frac12\phi''^{\,(0)}_{n_1,h_2-1}(\xi)\right] - \mu_2\left(1-\frac{\lambda_1}{\mu_1}\right)\left[\phi'^{\,(1)}_{n_1,h_2}(\xi)-\frac12\phi''^{\,(0)}_{n_1,h_2}(\xi)\right] \\
&+\omega\phi'^{\,(0)}_{n_1,h_2}(\xi),\\
\end{align*}
for $h_1=2,\dots,k_1$ and $h_2=k_1+2,\dots,k_1+k_2$. As we have done already a couple of times before, we can use the generating functions to easily sum all of these equations. Each equation contains the following three types of terms: $\phi'^{\,(0)}_{n_1,h}(\xi)$, $\phi'^{\,(1)}_{n_1,h}(\xi)$, and $\phi''^{\,(0)}_{n_1,h}(\xi)$. We denote the corresponding generating functions with $\PGF{Q}'^{\,(0)}_{h}(z, \xi)$, $\PGF{Q}'^{\,(1)}_{h}(z, \xi)$, and $\PGF{Q}''^{\,(0)}_{h}(z, \xi)$. After summing all equations, we obtain the following equation:
\begin{subequations}
\begin{align}
&\mu_2\frac{\lambda_1}{\mu_1}\sum_{h=1}^{k_1+k_2}\PGF{Q}'^{\,(1)}_{h}(1, \xi)-\mu_2\sum_{h_1=1}^{k_1}\PGF{Q}'^{\,(1)}_{h_1}(1, \xi) \label{finaleqn1}\\
&+\frac12\mu_2\sum_{h_2=k_1+1}^{k_1+k_2}\PGF{Q}''^{\,(0)}_{h_2}(1, \xi) +\frac12\mu_2\left(1-\frac{\lambda_1}{\mu_1}\right)\sum_{h=1}^{k_1+k_2}\PGF{Q}''^{\,(0)}_{h}(1, \xi)\label{finaleqn2}\\
&+\omega\sum_{h=1}^{k_1+k_2}\PGF{Q}'^{\,(0)}_{h}(1, \xi)\label{finaleqn3}\\
&=0.\nonumber
\end{align}
\end{subequations}

Note that the derivation of \eqref{finaleqn1} and \eqref{finaleqn2} follows closely the manner in which \eqref{solorder1} has been derived. The last term \eqref{finaleqn3} follows directly from collecting the $\omega$-terms.

Using the results obtained in appendices
the equation can be rewritten to
\begin{equation}
\frac{\lambda_1\mu_2^2}{\mu_1^2}P_0''(\xi)+\mu_2\left(1-\frac{\lambda_1}{\mu_1}\right)P''_0(\xi)+\omega P'_0(\xi)=0.\label{differentialeqnAlt}
\end{equation}
The proof that \eqref{finaleqn1} is the same as the first term in \eqref{differentialeqnAlt} can be found in Appendix \ref{proofOrder2}. The second and third term follow from \eqref{rhovacationsystem} and the fact that $\PGF{Q}^{\,(0)}_{h}(1, \xi) = \PGF{L}^{\,(0)}_{h}(1)P_0(\xi)$.
Equation \eqref{differentialeqnAlt} can be rewritten to \eqref{differentialeqn}, which concludes this proof.
\end{proof}
\end{proposition}

\subsection{The scaled number of customers in the critically loaded queue}

Now we can finally present the density of the scaled number of customers in $Q_2$, denoted by $P_0(\xi)$.
It is obtained by solving the differential equation \eqref{differentialeqn}:
\begin{equation}
P_0(\xi) = \eta \mathrm{e}^{-\eta \xi},
\end{equation}
with
\[
{\eta} = \omega\left[\mu_2+\frac{\lambda_1\mu_2(\mu_2-\mu_1)}{\mu_1^2}\right]^{-1}.
\]
We have used that $\int_0^\infty P_0(\xi)\textrm{d}\xi=1$ and that $\sum_{n_1=0}^\infty\sum_{n_2=0}^\infty\sum_{h=1}^{k_1+k_2}p(n_1,n_2,h)=1$. Without loss of generality, we may take $\omega = \mu_2$, which means that
\begin{equation}
\frac{1}{\eta} = 1-\frac{\lambda_1}{\mu_1}+\mu_2\frac{\lambda_1}{\mu_1^2}.\label{eta}
\end{equation}

We motivate this choice for $\omega$ by noting that we consider the scaled queue length $\delta N_2$ for $\delta\downarrow0$. By choosing $\omega=\mu_2$, and using \eqref{perturbation}, our scaling becomes equivalent to considering the scaled queue length $(1-\rho)N_2$, which is commonly used. Finally, by applying the multiclass distributional law of Bertsimas and Mourtzinou \cite{bertsimasmourtzinou} it directly follows that the scaled waiting time at $Q_2$ also follows an exponential distribution with parameter $\lambda_2 \eta$.

\subsection{Main result}

The analysis of the present section has the following immediate consequence for the joint (scaled) queue-length distribution in heavy traffic, which is the main result of this paper.

\paragraph{Main result:}
For $\lambda_1/k_1<\lambda_2/k_2$ and $\lambda_2=\mu_2\left((1-\rho_1)-\delta\right)$, we have:
\begin{equation}
\lim_{\delta\downarrow0}\P[N_1\leq n_1, \delta N_2\leq \xi] = \mathcal{L}(n_1)\left(1-\mathrm{e}^{-\eta \xi}\right),\label{mainresult}
\end{equation}
where $\mathcal{L}(\cdot)$ is the cumulative probability distribution of the queue length of a queueing system with multiple vacations satisfying Properties P\ref{vacationprop1}--P\ref{vacationproplast}, and $\eta$ is given by \eqref{eta}.

\section{Final remarks and suggestions for further research}\label{sect:conclusions}

\paragraph{Interpretation.}
The main result \eqref{mainresult} derived in the preceding section has the following intuitively appealing interpretation:
\begin{enumerate}
\item The number of customers in the stable queue has the same distribution as the number of customers in a $k$-limited vacation system with Erlang-$k_2$ distributed vacations.\label{propone}
\item The scaled number of customers in the critically loaded queue is exponentially distributed with parameter $\eta$.\label{proptwo}
\item The number of customers in the stable queue and the (scaled) number of customers in the critically loaded queue are independent.\label{propthree}
\end{enumerate}
Below we explain these properties heuristically.

Property \ref{propone} can be explained by the fact that if $Q_2$ is in heavy-traffic, then exactly $k_2$ customers are served at this queue during each
cycle. If we place an outside observer at $Q_1$, then, from his perspective, this queue behaves like a $k$-limited vacation model in heavy-traffic, where the vacation distribution is given by the convolutions of $k_2$ exponentially distributed service times distributions in $Q_2$.

For Property \ref{proptwo}, we note that the total workload in the system equals the amount of work in an $M/G/1$ queue in which the two customer classes are combined into one customer class with arrival rate $\lambda_1+\lambda_2$ and  hyperexponentially distributed service times, i.e., the service time is with probability $\lambda_i / (\lambda_1+\lambda_2)$ exponentially distributed with parameter $\mu_i$, $i=1,2$. Based on standard heavy-traffic results for the $M/G/1$ queue, this implies that the distribution of the scaled total workload converges to an exponential distribution with mean $\rho \E[R]$, where $R$ is a residual service time. For a hyperexponential distribution, we have
\[
\E[R]=\frac{\lambda_1/\mu_1^2+\lambda_2/\mu_2^2}{\rho},
\]
which implies that the total (scaled) asymptotic workload is exponentially distributed with parameter $\lambda_1/\mu_1^2+\lambda_2/\mu_2^2$.
In heavy traffic, almost all customers are located in $Q_2$ so the total number of customers at this queue is exponentially distributed with mean $\mu_2(\lambda_1/\mu_1^2+\lambda_2/\mu_2^2)$. Since $\lambda_2 \uparrow \mu_2(1-\lambda_1/\mu_1)$, the scaled number of customers in $Q_2$ is exponentially distributed with parameter $\eta$. Using the multiclass distributional law of Bertsimas and Mourtzinou \cite{bertsimasmourtzinou}, it can be shown that the scaled asymptotic waiting time of customers in $Q_2$ is exponentially distributed with parameter $\lambda_2\eta$.

Finally, Property \ref{propthree} follows from the time-scale separation in heavy-traffic which implies that the dynamics of the stable queue evolve at a much faster time scale than the dynamics of the critically loaded queue. Since the amount of ``memory'' of the stable queue asymptotically vanishes compared to that of the critically loaded queue, the number of customers in $Q_1$ and the scaled queue length of $Q_2$ become independent in the limit.

\paragraph{Two critically loaded queues.}
In the current paper we have analysed the heavy-traffic behaviour in case only $Q_2$ becomes critically loaded, i.e., when Assumption \eqref{stability2} is satisfied. The limiting regime in which both queues become saturated simultaneously $(\lambda_1/k_1=\lambda_2/k_2)$, shows fundamentally different system behaviour. That is, for general $\rho$ the waiting time at $Q_1$ is an (unknown) function of the visit time at $Q_2$ and $1/(1-\rho)$.  This implies that $\O(1-\rho)$ variations  in the visit time at $Q_2$ are relevant for the heavy-traffic behaviour at $Q_1$. More colloquially, it is not sufficient anymore to use a scaling that implies that exactly $k_2$ customers are served at $Q_2$ during each cycle, i.e., the probability that there are served less than $k_2$ customers cannot be neglected, when analysing the asymptotic behaviour of $Q_1$.

\paragraph{Further research.}

The analysis in this paper allows different kind of extensions. Firstly, one could consider phase-type interarrival-time or service-time distributions. The approach introduced in the present paper may be extended, without fundamentally changing the analysis, to such a system. Another extension could be the introduction of switch-over times whenever the server switches between queues. Such an extension requires more severe adaptations to the approach and the analysis, and is the topic of a forthcoming paper. Finally, we want to mention that the singular-perturbation technique can also be applied to derive the HT analysis of a system consisting of more than two (say $N$) queues, with one queue becoming critically loaded. Following the lines of the current paper, one can show that in this limiting regime the stable queues have the same joint queue-length distribution as in a polling model with $N-1$ queues and an extended switch-over time, whereas the scaled queue-length distribution of the critically loaded queue is again exponentially distributed. As such, the results of the present paper provide a theoretical basis for the transformation of large polling systems into smaller systems for approximation purposes, cf. LaPadula and Levy \cite{lapadulalevy96}.

\section*{Acknowledgement}

The authors would like to thank Sem Borst and Onno Boxma for interesting discussions and valuable comments on earlier drafts of this paper.

\appendix
\section*{Appendix}
\section{A vacation model with \boldmath$k$-limited service}\label{proofVacation}

In this appendix we study a queueing model with multiple vacations and $k$-limited service. The main goal is to find the PGF of the queue-length distribution, as to prove \eqref{lstL1}-\eqref{lstL3}. At the end of this appendix some additional properties of this queue-length distribution are given, which will be used in Section \ref{sect:equatingOdelta2} and in Appendix \ref{proofOrder2}.

The service times in this vacation model are exponentially distributed with parameter $\mu_1$, and the vacation length is Erlang($k_2$) distributed with parameter $\mu_2$. The service discipline is $k$-limited service, with at most $k_1$ customers being served during one visit period. Although the queue-length distribution for the case with generally distributed service and vacation times has been studied by Lee \cite{lee89}, we provide the proof here to keep the paper self-contained, but also because our state space is slightly different and we do not look at embedded epochs, yielding slightly different expressions than in \cite{lee89}.

The starting point is to obtain generating functions from the balance equations \eqref{balanceqn1order0}-\eqref{balanceqnlastorder0}. Multiplying Equation \eqref{balanceqnlastorder0} with $z^{n_1}$, summing over all $n_1=1,2,\dots$, and adding Equation \eqref{balanceqn2order0}, yields the following equation:
\begin{equation}
\LST{L}^{(0)}_{h_2}(z) = \LST{H}(z)\LST{L}^{(0)}_{h_2-1}(z),\label{lstSol1}
\end{equation}
for $h_2=k_1+2,\dots,k_1+k_2$, where $\LST{H}(z)$ is defined in \eqref{ghz}. The interpretation of \eqref{lstSol1} is that the number of customers in the system during a certain vacation stage is simply the number of customers present at the previous stage of the vacation, plus the arrivals during one (exponentially distributed) stage. Obviously, no customers leave the system during a vacation.

Multiplying Equation \eqref{balanceqn7order0} with $z^{n_1}$, summing over all $n_1=1,2,\dots$, and adding Equation \eqref{balanceqn1order0}, yields the following equation:
\begin{equation}
\LST{L}^{(0)}_{k_1+1}(z) = \LST{H}(z)\left(\frac{\mu_1}{\mu_2}\sum_{h=1}^{k_1-1}\pi^{(0)}_{1,h}+\pi^{(0)}_{0,k_1+k_2}
+ \frac{\mu_1}{\mu_2 z}\LST{L}^{(0)}_{k_1}(z)\right).
\end{equation}

Multiplying Equation \eqref{balanceqn6order0} with $z^{n_1+1}$, summing over all $n_1=1,2,\dots$, and adding Equation \eqref{balanceqn4order0} multiplied by $z$, yields the following equation:
\begin{equation}
\LST{L}^{(0)}_{h_1}(z) = \frac{\LST{G}(z)}{z}\left(\LST{L}^{(0)}_{h_1-1}(z)-\pi^{(0)}_{1,h_1-1}z\right),\label{lstSol3}
\end{equation}
for $h_1=2,3,\dots,k_1$, where $\LST{G}(z)$ is defined in \eqref{ghz}.

Multiplying Equation \eqref{balanceqn5order0} with $z^{n_1+1}$, summing over all $n_1=1,2,\dots$, and adding Equation \eqref{balanceqn3order0} multiplied by $z$, yields the following equation:
\begin{equation}
\LST{L}^{(0)}_{1}(z) = \LST{G}(z)\frac{\mu_2}{\mu_1}\left(\LST{L}^{(0)}_{k_1+k_2}(z)-\pi^{(0)}_{0,k_1+k_2}\right).\label{lstSol2}
\end{equation}

We now have $k_1+k_2$ equations, each of which expresses $\LST{L}^{(0)}_h(z)$ in terms of $\LST{L}^{(0)}_{h-1}(z)$ (and $\LST{L}^{(0)}_1(z)$ in terms of $\LST{L}^{(0)}_{k_1+k_2}(z)$).
Finally, we can solve these equations to find the expressions for $\LST{L}^{(0)}_{h}(z)$, for $h=1,2,\dots,k_1+k_2$. The results are given in \eqref{lstL1}-\eqref{lstL3}.

Note that there are still $k_1$ unknowns: $\pi^{(0)}_{1,h_1}$, for $h_1=1,\dots,k_1-1$, and $\pi^{(0)}_{0,k_1+k_2}$. These can be found using the roots of the denominator of \eqref{lstL1}. Rouch\'e's Theorem implies that the denominator has $k_1$ roots on and inside the unit circle. The requirement that \eqref{lstL1} should be regular inside the unit circle, implies that the numerator of \eqref{lstL1} should have these same roots \cite{winandsrouche}. Hence, we have a set of equations involving the roots and the numerator of \eqref{lstL1} to eliminate these $k_1$ unknowns.

\paragraph{Some additional properties.} In this paragraph we derive some results that are used throughout this paper, particularly in Section \ref{sect:equatingOdelta2} and in Appendix \ref{proofOrder2}. From a balancing argument, we know that the fraction of time that the system is in a vacation is $1-\rho$ (where, in this system, $\rho=\lambda_1/\mu_1$). Conversely, the fraction of time that the system is serving customers is $\rho$. Hence,
\begin{equation}
\sum_{h=1}^{k_1}\LST{L}^{(0)}_{h}(1) = \frac{\lambda_1}{\mu_1}, \qquad
\sum_{h=k_1+1}^{k_1+k_2}\LST{L}^{(0)}_{h}(1) = k_2\LST{L}^{(0)}_{k_1+k_2}(1) = 1-\frac{\lambda_1}{\mu_1}.\label{rhovacationsystem}
\end{equation}

Moreover, from \eqref{lstSol1} we know that $\LST{L}^{(0)}_{h_2}(1) = \LST{L}^{(0)}_{h_2-1}(1)$ for $h_2=k_1+2,\dots,k_1+k_2$. It follows that
\[
\LST{L}^{(0)}_{h_2}(1) =\frac{1}{k_2}\left(1-\frac{\lambda_1}{\mu_1}\right), \qquad h_2=k_1+1,\dots,k_1+k_2.
\]
From \eqref{lstSol3} and \eqref{lstSol2} we now have
\begin{equation}
\LST{L}^{(0)}_{h_1}(1) =\frac{\mu_2}{\mu_1}\left(\frac{1}{k_2}\left(1-\frac{\lambda_1}{\mu_1}\right)-\pi^{(0)}_{0,k_1+k_2}\right)-\sum_{i=1}^{h_1-1}\pi^{(0)}_{1,i}, \qquad h_1=1,\dots,k_1.\label{L0h1}
\end{equation}

The following relation for the unknowns $\pi^{(0)}_{0,k_1+k_2}$ and $\pi^{(0)}_{1,h_1}\quad (h_1=1,\dots,k_1-1)$ can be derived by combining all of these results:
\begin{equation}
k_1\pi^{(0)}_{0,k_1+k_2}+\frac{\mu_1}{\mu_2}\sum_{h=1}^{k_1-1}(k_1-h)\pi^{(0)}_{1,h} = \frac{1}{k_2}\left[k_1\left(1-\frac{\lambda_1}{\mu_1}\right)-k_2\frac{\lambda_1}{\mu_2}\right].\label{relationUnknowns0}
\end{equation}
This relation turns out to be crucial to derive many results in this paper, without having to know the exact expressions for all of the individual probabilities.

\begin{remark}
A balancing argument has been the starting point to derive all of the above properties. A rigorous way to derive these results, is by using L'H\^{o}pital's rule on \eqref{lstL1} to determine $\LST{L}^{(0)}_{k_1+k_2}(1)$, and subsequently deriving an expression for $\sum_{h=1}^{k_1+k_2}\LST{L}^{(0)}_h(1)$, which we know is equal to one.
\end{remark}

\section{The second-order perturbed balance equations}\label{proofOrder2}

The main goal of this appendix is to prove that \eqref{finaleqn1} can be written as the first term in \eqref{differentialeqnAlt}. If we rearrange the summations slightly, we can write the equation that we need to prove as follows:
\begin{equation}
\mu_2\frac{\lambda_1}{\mu_1}\sum_{h_2=k_1+1}^{k_1+k_2}\PGF{Q}'^{\,(1)}_{h_2}(1, \xi)-\mu_2\left(1-\frac{\lambda_1}{\mu_1}\right)\sum_{h_1=1}^{k_1}\PGF{Q}'^{\,(1)}_{h_1}(1, \xi)
=\frac{\lambda_1\mu_2^2}{\mu_1^2}P_0''(\xi).
\label{differentialeqnPart1}
\end{equation}

This equation should follow from Equations \eqref{balanceqn1order1}-\eqref{balanceqnlastorder1}. In order to prove it, we take the following steps:
\begin{enumerate}
\item First, we take the generating functions of Equations \eqref{balanceqn1order1}-\eqref{balanceqnlastorder1} to develop relations for $\PGF{Q}'^{\,(1)}_{h}(z, \xi) \quad (h=1,\dots,k_1+k_2)$.
\item The next step involves solving these equations to find an expression for $\PGF{Q}'^{\,(1)}_{k_1+k_2}(z, \xi)$.
\item Step~3 is to reformulate \eqref{differentialeqnPart1} in terms of $\PGF{Q}'^{\,(1)}_{k_1+k_2}(1, \xi)$. It turns out that in this stage all terms containing probabilities $\phi'^{\,(1)}_{n,h}(\xi)$ (or their generating functions) are eliminated.
\item The last step involves some more algebraic manipulations which eliminate all terms containing probabilities $\phi''^{\,(0)}_{n,h}(\xi)$ and, eventually, prove \eqref{differentialeqnPart1}.
\end{enumerate}

\paragraph{Step 1: Find relations for the generating functions.}

Multiplying Equation \eqref{balanceqnlastorder1} with $z^{n_1}$, summing over all $n_1=1,2,\dots$, and adding Equation \eqref{balanceqn2order1}, yields the following equation:
\begin{equation}
\LST{Q}^{(1)}_{h_2}(z, \xi) = \LST{H}(z)\left(\LST{Q}^{(1)}_{h_2-1}(z, \xi) + \LST{Q}'^{\,(0)}_{h_2-1}(z, \xi) -
\left(1-\frac{\lambda_1}{\mu_1}\right)\LST{Q}'^{\,(0)}_{h_2}(z, \xi)\right),
\label{eqnorder21}
\end{equation}
for $h_2=k_1+2,\dots,k_1+k_2$, where $\LST{H}(z)$ is defined in \eqref{ghz}.

Multiplying Equation \eqref{balanceqn7order1} with $z^{n_1}$, summing over all $n_1=1,2,\dots$, and adding Equation \eqref{balanceqn1order1}, yields the following equation:
\begin{align}
\LST{Q}^{(1)}_{k_1+1}(z, \xi) =& \LST{H}(z)\left(\frac{\mu_1}{\mu_2}\sum_{h=1}^{k_1-1}\phi^{(1)}_{1,h}(\xi) + \phi^{(1)}_{0,k_1+k_2}(\xi)+
\frac{\mu_1}{z\mu_2}\LST{Q}^{(1)}_{k_1}(z, \xi)+\phi'^{\,(0)}_{0,k_1+k_2}(\xi) \right.\nonumber\\
&\qquad\qquad\left.-\left(1-\frac{\lambda_1}{\mu_1}\right)\LST{Q}'^{\,(0)}_{k_1+1}(z, \xi)\right).
\label{eqnorder22}
\end{align}

Multiplying Equation \eqref{balanceqn6order1} with $z^{n_1+1}$, summing over all $n_1=1,2,\dots$, and adding Equation \eqref{balanceqn4order1} multiplied by $z$, yields the following equation:
\begin{equation}
\LST{Q}^{(1)}_{h_1}(z, \xi) = \frac{\LST{G}(z)}{z}\left(\LST{Q}^{(1)}_{h_1-1}(z, \xi)-\phi^{(1)}_{1,h_1-1}(\xi)z-\frac{z\mu_2}{\mu_1}\left(1-\frac{\lambda_1}{\mu_1}\right)\LST{Q}'^{\,(0)}_{h_1}(z, \xi)\right),
\label{eqnorder23}
\end{equation}
for $h_1=2,3,\dots,k_1$, where $\LST{G}(z)$ is defined in \eqref{ghz}.

Multiplying Equation \eqref{balanceqn5order1} with $z^{n_1+1}$, summing over all $n_1=1,2,\dots$, and adding Equation \eqref{balanceqn3order1} multiplied by $z$, yields the following equation:
\begin{align}
\LST{Q}^{(1)}_{1}(z, \xi) =& \frac{\mu_2}{\mu_1}\LST{G}(z)\left(\LST{Q}^{(1)}_{k_1+k_2}(z, \xi)-\phi^{(1)}_{0,k_1+k_2}(\xi)
+ \LST{Q}'^{\,(0)}_{k_1+k_2}(z, \xi)-\phi'^{(0)}_{0,k_1+k_2}(\xi)\phantom{\left(\frac{\lambda_1}{\mu_1}\right)}\right.\nonumber\\
&\qquad\qquad\left.-\left(1-\frac{\lambda_1}{\mu_1}\right)\LST{Q}'^{\,(0)}_{1}(z, \xi)\right).
\label{eqnorder24}
\end{align}

\paragraph{Step 2: Solve these relations and determine \boldmath $\PGF{Q}'^{\,(1)}_{k_1+k_2}(z, \xi)$.}

Solving Equations \eqref{eqnorder21}-\eqref{eqnorder24} requires a lot of straightforward, but tedious, computations. For reasons of compactness we will present some relevant intermediate results in this appendix, but leave the exact derivations to the reader. Firstly, one can use \eqref{eqnorder21}, combined with \eqref{lstSol1}, to express $\LST{Q}^{(1)}_{k_1+k_2}(z, \xi)$ in terms of $\LST{Q}^{(1)}_{k_1+1}(z, \xi)$ and $\LST{Q}'^{\,(0)}_{k_1+k_2}(z, \xi)$:
\begin{equation}
\LST{Q}^{(1)}_{k_1+k_2}(z, \xi) = \LST{H}(z)^{k_2-1}\LST{Q}^{(1)}_{k_1+1}(z, \xi)+(k_2-1)\left(1-\left(1-\frac{\lambda_1}{\mu_1}\right)\LST{H}(z)\right)\LST{Q}'^{\,(0)}_{k_1+k_2}(z, \xi).
\label{calc1}
\end{equation}
Second, we use \eqref{eqnorder22} to express $\LST{Q}^{(1)}_{k_1+1}(z, \xi)$ in terms of $\LST{Q}^{(1)}_{k_1}(z, \xi)$ and $\LST{Q}'^{\,(0)}_{k_1+1}(z, \xi)$. Subsequently, we use \eqref{eqnorder23} and \eqref{lstSol3} to express $\LST{Q}^{(1)}_{k_1}(z, \xi)$ in terms of $\LST{Q}^{(1)}_{1}(z, \xi)$ and $\LST{Q}'^{\,(0)}_{h}(z, \xi)$:
\begin{align}
\LST{Q}^{(1)}_{k_1}(z, \xi) &= \left(\frac{\LST{G}(z)}{z}\right)^{k_1-1}\LST{Q}^{(1)}_{1}(z, \xi)-\sum_{h=1}^{k_1-1}\phi^{(1)}_{1,h}(\xi)z\left(\frac{\LST{G}(z)}{z}\right)^{k_1-h} \nonumber\\
&- \frac{z\mu_2}{\mu_1}\left(1-\frac{\lambda_1}{\mu_1}\right)\sum_{h=2}^{k_1}\left(\frac{\LST{G}(z)}{z}\right)^{k_1-h+1}\LST{Q}'^{\,(0)}_{h}(z, \xi).
\label{calc2}
\end{align}

Finally, we use \eqref{eqnorder24} to express $\LST{Q}^{(1)}_{1}(z, \xi)$ in terms of $\LST{Q}^{(1)}_{k_1+k_2}(z, \xi)$ again. After some rearrangement of the terms, this leads to the following expression for $\LST{Q}^{(1)}_{k_1+k_2}(z, \xi)$:
\begin{equation}
\LST{Q}^{(1)}_{k_1+k_2}(z, \xi) = \frac{A^{(1)}(z,\xi)+A'^{\,(0)}(z,\xi)}{D(z)},\label{lstQ2}
\end{equation}
where
\begin{align}
A^{(1)}(z,\xi) &= \displaystyle
\LST{H}(z)^{k_2}\left[\phi^{(1)}_{0,k_1+k_2}(\xi) \left(1- \left({\LST{G}(z)}/{z}\right)^{k_1}\right) +\frac{\mu_1}{\mu_2}\sum_{h=1}^{k_1-1}\phi_{1,h}^{(1)}(\xi)\left(1-\left({\LST{G}(z)}/{z}\right)^{k_1-h} \right) \right],
\label{A1}\\
A'^{\,(0)}(z,\xi) &= k_2 \LST{Q}'^{\,(0)}_{k_1+k_2}(z, \xi) - \LST{H}(z)^{k_2}\frac{\mu_1}{\mu_2}\sum_{h=1}^{k_1-1}\phi'^{\,(0)}_{1,h}(\xi)\left(1-\left({\LST{G}(z)}/{z}\right)^{k_1-h} \right) \nonumber\\
&-\left(1-\frac{\lambda_1}{\mu_2}\right)\left[\LST{H}(z)^{k_2}\sum_{h=1}^{k_1}\LST{Q}'^{\,(0)}_{h}(z, \xi)\left({\LST{G}(z)}/{z}\right)^{k_1-h+1}   + k_2\LST{H}(z)\LST{Q}'^{\,(0)}_{k_1+k_2}(z, \xi)\right],
\label{A0}\\
D(z)&=1-\big(\LST{G}(z)/z\big)^{k_1}\LST{H}(z)^{k_2}\label{D}.
\end{align}

Note that $A^{(1)}(z,\xi)$ only contains the probabilities $\phi^{(1)}_{n,h}(\xi)$. All probabilities $\phi'^{\,(0)}_{0,k_1+k_2}(\xi)$, and their generating functions, are contained in $A'^{\,(0)}(z,\xi)$. Also note that $\frac{A^{(1)}(z,\xi)}{D(z)}$ is exactly the same expression as \eqref{lstL1}, but with constants $\pi_{n,h}^{(0)}$ replaced by $\phi_{n,h}^{(1)}$. The reason is that, if one would ignore the probabilities $\phi'^{\,(0)}_{n,h}(\xi)$ in the balance equations \eqref{balanceqn1order1}-\eqref{balanceqnlastorder1}, the system is completely equivalent to the system \eqref{balanceqn1order0}-\eqref{balanceqnlastorder0}, which corresponds to the vacation system studied in Appendix \ref{proofVacation}.

\paragraph{Step 3: reformulate the original problem in terms of \boldmath$\PGF{Q}'^{\,(1)}_{k_1+k_2}(1, \xi)$.}
In order to solve \eqref{differentialeqnPart1}, we need to determine $\sum_{h_1=1}^{k_1}\PGF{Q}'^{\,(1)}_{h_1}(1, \xi)$ and $\sum_{h_1=k_1+1}^{k_1+k_2}\PGF{Q}'^{\,(1)}_{h_1}(1, \xi)$. After substituting $z=1$ in Equations \eqref{eqnorder21}-\eqref{eqnorder24}, we can express these sums in terms of $\LST{Q}'^{\,(1)}_{k_1+k_2}(1, \xi)$ and $\LST{Q}''^{\,(0)}_{h}(1, \xi)$:
\begin{align}
\sum_{h_1=1}^{k_1}\PGF{Q}'^{\,(1)}_{h_1}(1, \xi) &= k_1\frac{\mu_2}{\mu_1}\left(\LST{Q}'^{\,(1)}_{k_1+k_2}(1, \xi)-\phi'^{\,(1)}_{0,k_1+k_2}\right)
-\sum_{h=1}^{k_1-1}(k_1-h)\phi'^{\,(1)}_{1,h}\nonumber\\
&+ k_1 \PGF{Q}''^{\,(0)}_{1}(1, \xi)-\frac{\mu_2}{\mu_1}\left(1-\frac{\lambda_1}{\mu_1}\right)\sum_{h=1}^{k_1}(k_1-h+1)\PGF{Q}''^{\,(0)}_{h}(1, \xi),
\label{sum1}\\
\sum_{h_1=k_1+1}^{k_1+k_2}\PGF{Q}'^{\,(1)}_{h_1}(1, \xi) &= k_2\PGF{Q}'^{\,(1)}_{k_1+k_2}(1, \xi) - \frac{\lambda_1}{\mu_1}\left(\frac{k_2-1}{2}\right)\left(1-\frac{\lambda_1}{\mu_1}\right)P''_0(\xi).
\label{sum2}
\end{align}
Since $\PGF{Q}''^{\,(0)}_{h}(1, \xi)$ for $h=1,\dots,k_1$ can be determined directly using \eqref{L0h1}, it only remains to determine $\PGF{Q}'^{\,(1)}_{k_1+k_2}(1, \xi)$. Fortunately, according to the following lemma we can focus on the part $\lim_{z\rightarrow1}\frac{A'^{\,(0)}(z,\xi)}{D(z)}$ only.

\begin{lemma}\label{lemmaOrder1termscancel}
The probabilities $\phi'^{\,(1)}_{0,k_1+k_2}$ and $\phi'^{\,(1)}_{0,h}$ $(h=1,\dots,k_1-1)$ in the left hand side of Equation \eqref{differentialeqnPart1} cancel out. Using \eqref{sum1} and \eqref{sum2}, we can express this statement in a more formal presentation:
\begin{align}
&\lim_{z\rightarrow1}\mu_2\frac{\lambda_1}{\mu_1}\left(k_2\frac{A^{(1)}(z,\xi)}{D(z)}\right) - \mu_2\left(1-\frac{\lambda_1}{\mu_1}\right)\left(k_1\frac{\mu_2}{\mu_1}\left(\frac{A^{(1)}(z,\xi)}{D(z)}-\phi^{(1)}_{0,k_1+k_2}\right)
-\sum_{h=1}^{k_1-1}(k_1-h)\phi^{(1)}_{1,h}\right)=0.
\label{order1termscancel}
\end{align}
\begin{proof}
Using L'H\^{o}pital's rule it can be shown that
\begin{equation}
\lim_{z\rightarrow1}\frac{A^{(1)}(z,\xi)}{D(z)} = \frac{\left(1-\frac{\lambda_1}{\mu_1}\right)C^{(1)}}{k_1\left(1-\frac{\lambda_1}{\mu_1}\right)-k_2\frac{\lambda_1}{\mu_2}},
\label{limit1}
\end{equation}
where
\[
C^{(1)} = k_1\phi^{(1)}_{0,k_1+k_2}+\frac{\mu_1}{\mu_2}\sum_{h=1}^{k_1-1}(k_1-h)\phi^{(1)}_{1,h}.
\]
After substituting \eqref{limit1} in \eqref{order1termscancel}, it is easily shown that the left hand side of \eqref{order1termscancel} indeed equals zero.
\end{proof}
\end{lemma}

\paragraph{Step 4: solve the original problem.}
We are almost ready to solve Equation~\eqref{differentialeqnPart1}, but we compute two helpful intermediate results first.
\begin{lemma}\label{lemmaY}
Define $Y := \lim_{z\rightarrow1}\frac{\partial}{\partial z}\LST{Q}'^{\,(0)}_{k_1+k_2}(z, \xi)$. This can be written as
\begin{align}
Y =&
\left[{2\frac{\mu_2}{\mu_1}\left(k_1\left(1-\frac{\lambda_1}{\mu_1}\right)-k_2\frac{\lambda_1}{\mu_2}\right)}\right]^{-1}
\times\left[
\left(1-\frac{\lambda_1}{\mu_1}\right)^2\sum_{h=1}^{k_1-1}h(k_1-h)\phi'^{\,(0)}_{1,h} + \right.\nonumber\\
&\qquad\left.
\left(\frac{\lambda_1}{\mu_1}\left(1-\frac{\lambda_1}{\mu_1}\right)^2(k_1+1)-\frac{\lambda_1^2}{\mu_1\mu_2}\left(1-\frac{\lambda_1}{\mu_1}\right)(k_2-1)
+2\left(\frac{\lambda_1}{\mu_1}\right)^2\right)P_0'(\xi)
\right].
\label{dL0k1k2}
\end{align}
\begin{proof}
Equation \eqref{dL0k1k2} follows after differentiating \eqref{lstL1} with respect to $z$, and subsequently applying L'H\^{o}\-pital's rule twice. Obviously, some basic algebraic manipulations are required to obtain the presentation in \eqref{dL0k1k2}.
\end{proof}
\end{lemma}

\begin{lemma}\label{lemmaX}
Define $X := \lim_{z\rightarrow1}\frac{A'^{\,(0)}(z,\xi)}{D(z)}$. This can be written as
\begin{align}
X&=-\frac{\mu_2}{\mu_1}Y +\left[k_1\left(1-\frac{\lambda_1}{\mu_1}\right)-k_2\frac{\lambda_1}{\mu_2}\right]^{-1}\times  \nonumber\\
&\left\{\phi'^{\,(0)}_{0,k_1+k_2}\left(k_1\left(1-\frac{\lambda_1}{\mu_1}\right)\left(k_2\frac{\lambda_1}{\mu_1}+\left(\frac{\lambda_1}{\mu_1}-k_1\left(1-\frac{\lambda_1}{\mu_1}\right)\right)\frac{\mu_2}{\mu_1}\right)\right)
\right.\nonumber\\
&\ +\left(1-\frac{\lambda_1}{\mu_1}\right)\left(k_2\frac{\lambda_1}{\mu_2}+\frac{\lambda_1}{\mu_1}-\frac{\mu_1}{\mu_2}-k_1\left(1-\frac{\lambda_1}{\mu_1}\right)\right)\sum_{h=1}^{k_1-1}(k_1-h)\phi'^{\,(0)}_{1,h}\nonumber\\
&\ +\left(1-\frac{\lambda_1}{\mu_1}\right)^2\,\sum_{h=1}^{k_1-1}h(k_1-h)\phi'^{\,(0)}_{1,h}\nonumber\\
&\ \left.-\left(1-\frac{\lambda_1}{\mu_1}\right)^2\,\left(\frac{k_1}{k_2}\left(\frac{\lambda_1}{\mu_1}-k_1\left(1-\frac{\lambda_1}{\mu_1}\right)\right)\frac{\mu_2}{\mu_1}+\frac{\lambda_1}{\mu_1}\left(\frac{\mu_1}{\mu_2}+k_1\right)\right)P_0'(\xi)\right\}.
\end{align}
\begin{proof}
This equation follows from applying L'H\^opital's rule to $\frac{A'^{\,(0)}(z,\xi)}{D(z)}$ and, hence, differentiating \eqref{A0} and \eqref{D} with respect to $z$. Substitution of $z=1$ gives the desired result after, again, many algebraic manipulations.
\end{proof}
\end{lemma}

Finally, we are ready to present the main result of this appendix, which is the proof of Equation \eqref{differentialeqnPart1}. Using Lemma \ref{lemmaOrder1termscancel}, and Equations \eqref{sum1} and \eqref{sum2}, we can write the left hand side of Equation \eqref{differentialeqnPart1} as
\begin{align*}
&\mu_2\frac{\lambda_1}{\mu_1}\left(k_2 X'-\frac{\lambda_1}{\mu_1}\left(\frac{k_2-1}{2}\right)\left(1-\frac{\lambda_1}{\mu_1}\right)P''_0(\xi)\right)\nonumber\\
&-\mu_2\left(1-\frac{\lambda_1}{\mu_1}\right)\left[k_1\frac{\mu_2}{\mu_1}X'+ k_1 \PGF{Q}''^{\,(0)}_{1}(1, \xi)-\frac{\mu_2}{\mu_1}\left(1-\frac{\lambda_1}{\mu_1}\right)\sum_{h=1}^{k_1}(k_1-h+1)\PGF{Q}''^{\,(0)}_{h}(1, \xi)\right],
\end{align*}
where $X'$ is the derivative of $X$ with respect to $\xi$. Using Lemma~\ref{lemmaX} and Lemma~\ref{lemmaY} (and Equation~\eqref{L0h1} to determine $\PGF{Q}''^{\,(0)}_{h}(1, \xi)$ for $h=1,\dots,k_1$) we can show that the above expression reduces to $\frac{\lambda_1\mu_2^2}{\mu_1^2}P_0''(\xi)$.

\expandafter\ifx\csname urlstyle\endcsname\relax
  \providecommand{\doi}[1]{DOI: #1}\else
  \providecommand{\doi}{DOI: \begingroup \urlstyle{rm}\Url}\fi

\bibliographystyle{abbrvnat}

\begin{thebibliography}{28}
\providecommand{\natexlab}[1]{#1}
\providecommand{\url}[1]{\texttt{#1}}
\expandafter\ifx\csname urlstyle\endcsname\relax
  \providecommand{\doi}[1]{doi: #1}\else
  \providecommand{\doi}{doi: \begingroup \urlstyle{rm}\Url}\fi

\bibitem[Adan et~al.(2006)Adan, van Leeuwaarden, and Winands]{winandsrouche}
I.~J. B.~F. Adan, J.~S.~H. van Leeuwaarden, and E.~M.~M. Winands.
\newblock On the application of {Rouch\'e's} theorem in queueing theory.
\newblock \emph{Operations Research Letters}, 34\penalty0 (3):\penalty0
  355--360, 2006.

\bibitem[Bertsimas and Mourtzinou(1997)]{bertsimasmourtzinou}
D.~Bertsimas and G.~Mourtzinou.
\newblock Multiclass queueing systems in heavy traffic: {An} asymptotic
  approach based on distributional and conservation laws.
\newblock \emph{Operations Research}, 45\penalty0 (3):\penalty0 470--487, 1997.

\bibitem[Blanc(1992)]{blanc92}
J.~P.~C. Blanc.
\newblock An algorithmic solution of polling models with limited service
  disciplines.
\newblock \emph{IEEE Transactions of Communications}, 40\penalty0 (7):\penalty0
  1152--1155, 1992.

\bibitem[Boon et~al.(2012)Boon, Adan, Winands, and Down]{boontrafficlights2012}
M.~A.~A. Boon, I.~J. B.~F. Adan, E.~M.~M. Winands, and D.~G. Down.
\newblock Delays at signalised intersections with exhaustive traffic control.
\newblock \emph{Probability in the Engineering and Informational Sciences},
  26\penalty0 (3):\penalty0 337--373, 2012.

\bibitem[Borst et~al.(1995)Borst, Boxma, and Levy]{borst}
S.~C. Borst, O.~J. Boxma, and H.~Levy.
\newblock The use of service limits for efficient operation of multistation
  single-medium communication systems.
\newblock \emph{IEEE/ACM Transactions on Networking}, 3\penalty0 (5):\penalty0
  602--612, 1995.

\bibitem[Boxma(1985)]{boxma1}
O.~J. Boxma.
\newblock Two symmetric queues with alternating service and switching times.
\newblock In E.~Gelenbe, editor, \emph{Performance '84}, pages 409--431.
  North-Holland, Amsterdam, 1985.

\bibitem[Boxma and Groenendijk(1988)]{boxma2}
O.~J. Boxma and W.~P. Groenendijk.
\newblock Two queues with alternating service and switching times.
\newblock In O.~J. Boxma and R.~Syski, editors, \emph{Queueing Theory and its
  Applications - Liber Amicorum for J.W. Cohen}, pages 261--282. North-Holland,
  Amsterdam, 1988.

\bibitem[Charzinski et~al.(1994)Charzinski, Renger, and Tangemann]{charzinski}
J.~Charzinski, T.~Renger, and M.~Tangemann.
\newblock Simulative comparison of the waiting time distributions in cyclic
  polling systems with different service strategies.
\newblock In \emph{Proceedings of the 14th International Teletraffic Congress,
  Antibes Juan-les-Pins}, pages 719--728, 1994.

\bibitem[{Coffman, Jr.} et~al.(1995){Coffman, Jr.}, Puhalskii, and
  Reiman]{coffman95}
E.~G. {Coffman, Jr.}, A.~A. Puhalskii, and M.~I. Reiman.
\newblock Polling systems with zero switchover times: A heavy-traffic averaging
  principle.
\newblock \emph{The Annals of Applied Probability}, 5\penalty0 (3):\penalty0
  681--719, 1995.

\bibitem[{Coffman, Jr.} et~al.(1998){Coffman, Jr.}, Puhalskii, and
  Reiman]{coffman98}
E.~G. {Coffman, Jr.}, A.~A. Puhalskii, and M.~I. Reiman.
\newblock Polling systems in heavy-traffic: A {Bessel} process limit.
\newblock \emph{Mathematics of Operations Research}, 23:\penalty0 257--304,
  1998.

\bibitem[Cohen and Boxma(1981)]{cohen1}
J.~W. Cohen and O.~J. Boxma.
\newblock The {$M/G/1$} queue with alternating service formulated as a
  {Riemann-Hilbert} problem.
\newblock In F.~J. Kylstra, editor, \emph{Performance '81}, pages 181--189,
  1981.

\bibitem[Down(1998)]{down}
D.~G. Down.
\newblock On the stability of polling models with multiple servers.
\newblock \emph{Journal of Applied Probability}, 35:\penalty0 925--935, 1998.

\bibitem[Eisenberg(1979)]{eisenberg1}
M.~Eisenberg.
\newblock Two queues with alternating service.
\newblock \emph{SIAM Journal on Applied Mathematics}, 36\penalty0 (2):\penalty0
  287--303, 1979.

\bibitem[Fricker and Ja\"ibi(1994)]{fricker}
C.~Fricker and M.~R. Ja\"ibi.
\newblock Monotonicity and stability of periodic polling models.
\newblock \emph{Queueing Systems}, 15:\penalty0 211--238, 1994.

\bibitem[Fuhrmann(1981)]{fuhrmann}
S.~W. Fuhrmann.
\newblock Performance analysis of a class of cyclic schedules.
\newblock Technical Memorandum 81-59531-1, Bell Laboratories, 1981.

\bibitem[Groenendijk(1990)]{Groenendijk1}
W.~P. Groenendijk.
\newblock \emph{Conservation Laws in Polling Systems}.
\newblock PhD thesis, University of Utrecht, 1990.

\bibitem[Ibe(1990)]{ibe}
O.~C. Ibe.
\newblock Analysis of polling systems with mixed service disciplines.
\newblock \emph{Stochastic Models}, 6:\penalty0 667--689, 1990.

\bibitem[Keilson and Servi(1990)]{keilsonservi90}
J.~Keilson and L.~D. Servi.
\newblock The distributional form of {Little's Law} and the {Fuhrmann-Cooper}
  decomposition.
\newblock \emph{Operations Research Letters}, 9\penalty0 (4):\penalty0
  239--247, 1990.

\bibitem[Knessl and Tier(1995)]{knessltiersurvey}
C.~Knessl and C.~Tier.
\newblock Applications of singular perturbation methods in queueing.
\newblock In \emph{Advances in Queueing Theory, Methods, and Open Problems},
  Probability and Stochastics Series, pages 311--336. CRC Press, 1995.

\bibitem[LaPadula and Levy(1996)]{lapadulalevy96}
C.~A. LaPadula and H.~Levy.
\newblock Customer delay in very large multi-queue single-server systems.
\newblock \emph{Performance Evaluation}, 26\penalty0 (3):\penalty0 201--218,
  1996.

\bibitem[Lee(1996)]{lee}
D.-S. Lee.
\newblock A two-queue model with exhaustive and limited service disciplines.
\newblock \emph{Stochastic Models}, 12\penalty0 (2):\penalty0 285--305, 1996.

\bibitem[Lee(1989)]{lee89}
T.~T. Lee.
\newblock {$M/G/1/N$} queue with vacation time and limited service discipline.
\newblock \emph{Performance Evaluation}, 9\penalty0 (3):\penalty0 181--190,
  1989.

\bibitem[Morrison and Borst(2010)]{morrisonborst2010}
J.~A. Morrison and S.~C. Borst.
\newblock Interacting queues in heavy traffic.
\newblock \emph{Queueing Systems}, 65\penalty0 (2):\penalty0 135--156, 2010.

\bibitem[Ozawa(1990)]{ozawa1}
T.~Ozawa.
\newblock Alternating service queues with mixed exhaustive and {$K$}-limited
  services.
\newblock \emph{Performance Evaluation}, 11:\penalty0 165--175, 1990.

\bibitem[Ozawa(1997)]{ozawa2}
T.~Ozawa.
\newblock Waiting time distribution in a two-queue model with mixed exhaustive
  and gated-type {$K$}-limited services.
\newblock In \emph{Proceedings of International Conference on the Performance
  and Management of Complex Communication Networks}, pages 231--250, 1997.

\bibitem[Resing(1993)]{resing}
J.~A.~C. Resing.
\newblock Polling systems and multitype branching processes.
\newblock \emph{Queueing Systems}, 13:\penalty0 409--426, 1993.

\bibitem[Van~der Mei(2007)]{RvdM_QUESTA}
R.~D. Van~der Mei.
\newblock Towards a unifying theory on branching-type polling models in heavy
  traffic.
\newblock \emph{Queueing Systems}, 57:\penalty0 29--46, 2007.

\bibitem[Winands et~al.(2009)Winands, Adan, van Houtum, and Down]{winands}
E.~M.~M. Winands, I.~J. B.~F. Adan, G.~J. van Houtum, and D.~G. Down.
\newblock A state-dependent polling model with $k$-limited service.
\newblock \emph{Probability in the Engineering and Informational Sciences},
  23\penalty0 (2):\penalty0 385--408, 2009.

\end{thebibliography}

\end{document}